\documentclass{kms-j}

\newcommand{\publname}{Preprint arXiv}
\issueinfo{}
  {}
  {}
  {}
\pagespan{1}{}
\copyrightinfo{}
  {Korean Mathematical Society}

\usepackage{graphicx}
\usepackage{epstopdf}              
\usepackage[style=numeric, backend=bibtex]{biblatex}
\allowdisplaybreaks

\theoremstyle{plain}
\newtheorem{theorem}{Theorem}[section]

\newtheorem{lemma}[theorem]{Lemma}
\newtheorem{corollary}[theorem]{Corollary}

\theoremstyle{definition}

\theoremstyle{remark}
\newtheorem{remark}[theorem]{Remark}

\begin{document}

\title[Bohr's Phenomenon]{Bohr's Phenomenon for Some Univalent Harmonic Functions}

\author[C. Singla]{Chinu Singla}
\address{Chinu Singla \\ Department of Mathematics \\ Sant Longowal Institute of Engineering and Technology\\ Longowal 148106, Punjab\\ India.}
\email{chinusingla204@gmail.com}

\author[S. Gupta]{Sushma Gupta}
\address{Sushma Gupta \\ Department of Mathematics \\ Sant Longowal Institute of Engineering and Technology\\ Longowal 148106, Punjab\\ India.}
\email{sushmagupta1@yahoo.com}

\author[S. Singh]{Sukhjit Singh}
\address{Sukhjit Singh \\ Department of Mathematics \\ Sant Longowal Institute of Engineering and Technology\\ Longowal 148106, Punjab\\ India.}
\email{sukhjit\_d@yahoo.com}
\thanks{Work of first author is supported by UGC - CSIR, New Delhi in the form of Senior Research Fellowship vide Award Letter Number - 2061540842.}

\subjclass{Primary 30A10, 30C45}
\keywords{Bohr radius, harmonic univalent functions, convex in one direction}

\begin{abstract}
In 1914 Bohr proved that there is an $r_0 \in(0,1)$ such that if a power series $\sum_{m=0}^\infty c_m z^m$  is convergent in the open unit disc and $|\sum_{m=0}^\infty c_m z^m|<1$ then, $\sum_{m=0}^\infty |c_m z^m|<1$ for $|z|<r_0$. The largest value of such $r_0$ is called the Bohr radius. In this article, we find Bohr radius for some univalent harmonic mappings having different dilatations and in addition, also compute Bohr radius for the functions convex in one direction.
\end{abstract}

\maketitle

\section{Introduction}
The Bohr inequality, first introduced in 1914 by Harald Bohr in his seminal work \cite{bohr1914theorem} and subsequently improved independently by M. Riesz, I. Shur and F. Wiener, essentially states that if $f(z)=\sum_{m=0}^\infty a_m z^m$  is an analytic function in the open unit disc $\mathbb{D}=\{z \in \mathbb{C}:|z|<1\}$  and $|f(z)|<1$ for all  $z\in {\mathbb{D}}$, then
\begin{equation}\label{eq:1}
\sum_{m=0}^\infty |a_m| r^m \leq 1
\end{equation}
for all  $z\in {\mathbb{D}}$ with $|z|=r\leq r_0=1/3$ and $1/3$ is the largest possible value of $r_0$, called the Bohr radius. The inequalities of the type \eqref{eq:1} have become famous by the name {\it{ Bohr inequalities}} and the problems of finding the largest possible values of $r_0$ in different setups are now a days called  the {\it Bohr radius problems}. For a glimpse of the ongoing current research in this area we refer the reader to some recent articles, e.g. \cite{bhowmik2019bohr, bhowmik2018bohr, kayumov2018bohr, muhanna2010bohr, muhanna2017bohr} and references therein.\\
In 2010, Abu Muhanna \cite{muhanna2010bohr} investigated some Bohr radius problems using the concept of subordination. For two analytic functions $f$ and $g$ in $\mathbb{D},$ $g$ is said to be subordinate to $f$( written, $g\prec f)$ if there exists a function $\psi$ analytic in $\mathbb{D}$ with $\psi(0)=0$ and $|\psi(z)|<1$ such that $g=fo\psi.$ In particular, when $f$ is univalent, then $g\prec f$ is equivalent to $g(0)=f(0)$ and $g(\mathbb{D})\subset f(\mathbb{D})$. We shall denote by $S(f),$ the class of  all functions $g$ subordinate to a fixed function $f.$ A class of analytic (harmonic) functions in the unit disc $\mathbb{D}$ is said to possess classical Bohr's phenomenon if an inequalty of the type \eqref{eq:1} is satisfied in $|z|<r_0,$ for some $r_0,0<r_0\leq 1$. It is known (see \cite{muhanna2010bohr}) that not all classes of functions have classical Bohr's phenomenon. So, Abu Muhanna \cite{muhanna2010bohr} reformulated classical Bohr's phenomenon and proved the following result:\\
\begin{theorem}

If $f(z)=\sum_{m=0}^\infty a_m z^m$ is a univalent function and $g(z)=\sum_{m=0}^\infty b_m z^m \in S(f),$ then
\begin{equation} \label{eq:2}
\sum_{m=1}^\infty |b_m| r^m\leq d(f(0),\partial f(\mathbb{D}))          
\end{equation}
for all $|z|=r\leq r_0=3-\sqrt{8}=0.17157..,$ where $d(f(0),\partial f(\mathbb{D}))$ is the Euclidean distance between $f(0)$ and $\partial f(\mathbb{D}),$ the boundary of $f(\mathbb{D}).$ The value of $r_0$ is sharp for $f(z)=z/(1-z)^2 $, the Koebe function. Further, if $f$ is convex univalent in $\mathbb{D}$, then $r_0=1/3.$ 
\end{theorem}
In the recent years, a number of research articles (for example see \cite{bhowmik2019bohr}, \cite{liu2020bohr}, \cite{liu2019bohr}) are published and many hidden facts of this subject are brought into broad daylight. In particular, Bhowmik and Das \cite{bhowmik2019bohr} successfully extended the Bohr inequalities of type \eqref{eq:2} for certain harmonic functions. A complex valued function $f(z)=u(x,y)+iv(x,y)$ of $z=x+iy \in \mathbb{D}$ is said to be harmonic if both $u(x,y)$ and $v(x,y)$ are  real harmonic in $\mathbb{D}$. It is known that such an $f$ can be uniquely represented as $f=h+\overline{g},$ where $h$ and $g$ are analytic functions in $\mathbb{D}$ with $f(0)=h(0)$. It immediately follows from this representation that  $f$ is locally univalent and sense preserving whenever its Jacobian $J_f$, defined by $J_f(z):=|h'(z)|-|g'(z)|$, satisfies $J_f(z)>0$ for all $z\in \mathbb{D};$ or equivalently if $h'\ne 0$ in $\mathbb{D}$ and the (second complex)  dilatation $w_f$ of $f$, defined by $w_f(z)=g'(z)/h'(z)$, satisfies the condition $|w_f(z)|<1$ in $\mathbb{D}$. A  harmonic function $f=h+\overline{g}$ defined in $\mathbb{D}$ is said to be K-quasiconformal if its dilatation $w_f$ satisfies $|w_f|\leq k, k=(K-1)/(K+1)\in [0,1).$ In view of the work of Schaubroeck in \cite{schaubroeck2000subordination}, aforesaid definitions and notations for subordination of analytic functions can be extended to harmonic functions without any change. This lead Bhowmik and Das \cite{bhowmik2019bohr} to extend Theorem 1.1 as under:
\begin{theorem}
Let $f(z)=h(z)+\overline{g(z)}=\sum_{m=0}^{\infty} a_m z^m + \overline{\sum_{m=1}^{\infty} b_m z^m}$ be a sense preserving K-quasiconformal harmonic mapping defined in $\mathbb{D}$ such that $h$ is univalent in $\mathbb{D},$ and let $f_1(z)=h_1(z)+\overline{g_1(z)}=\sum_{m=0}^{\infty} c_m z^m + \overline{\sum_{m=0}^{\infty} d_m z^m} \in S(f).$ Then
\begin{equation}\label{eq:3}
\sum_{m=1}^\infty |c_m| r^m+\sum_{m=1}^\infty |d_m| r^m \leq d(h(0),\partial h(\mathbb{D})) 
\end{equation}
for $|z|=r\leq r_0=(5K+1-\sqrt{8K(3K+1)}/(K+1).$ This result is sharp for the function $p(z)=z/(1-z)^2+k \overline{z/(1-z)^2},$ where $k=(K-1)/(K+1).$ Moreover, if we take $h$ to be convex univalent then the inequality in (\ref{eq:3}) holds for $|z|=r\leq r_0=(K+1)/(5K+1).$ This result is again sharp for the function $q(z)=z/(1-z)+k \overline{z/(1-z)}.$
\end{theorem}
In this article, our aim is to establish the Bohr's phenomenon and compute Bohr radius for some subclasses of univalent harmonic functions. We also propose to improvise Theorem 1.1 and 1.2  stated above. \\

We close this section by setting certain notations for subsequent use in this paper. We denote by $S_H,$ the class of univalent harmonic functions $f$ normalized by the conditions $f(0)=0$ and $f_z(0)=1.$ In addition, if $f_{\overline{z}}(0)=0$ also, then the class is denoted by $S_H^0.$ Further, $K_H^0$ is the usual subclass of $S_H^0$ consisting of convex functions. A domain $\Omega$ is said to be convex in the direction $\theta, 0\le \theta<\pi$, if the intersection of the straight line through the origin and the point $e^{i\theta}$ in the complex plane is connected or empty. A function $f$ mapping the open unit disc $\mathbb{D}$ onto such a domain is called \textit{convex in direction $\theta$.}\\
\section{Main Results}
We begin this section by stating following lemma which immediately follows from the work of Bhowmik and Das \cite{bhowmik2018bohr}.
\begin{lemma}
Let $f(z)=\sum_{m=0}^{\infty} a_m z^m$ and $g(z)=\sum_{m=0}^{\infty} b_m z^m$ be two analytic functions in $\mathbb{D}$ and $g\prec f.$ Then
$$ \sum_{m=0}^{\infty} |b_m |r^m \leq \sum_{m=0}^{\infty} |a_m| r^m$$
for $|z|=r\leq 1/3.$
\end{lemma}
Using this lemma, we now improvise Theorem 1.1 by taking univalent analytic function in $\mathbb{D}$ as $f(z)=z+\sum_{m=2}^\infty a_m z^m.$ Making use of well known De Brange's theorem: $|a_m|\leq m, m=2,3,...,$ and after some simple calculations, we easily get: 
\begin{theorem}
If $f(z)=z+\sum_{m=2}^\infty a_m z^m$ is a univalent analytic function in $\mathbb{D}$ and $g(z)=\sum_{m=1}^\infty b_m z^m \in S(f),$ then 
\begin{equation}\label{eq:2.1}
\sum_{m=1}^\infty |b_m| r^m\leq 1
\end{equation} 
for all $|z|=r\leq 1/3.$
\end{theorem}
 
In a similar manner, we restate Theorem 1.2 as under;

\begin{theorem}
Let $f(z)=h(z)+\overline{g(z)}=z+\sum_{m=2}^{\infty} a_m z^m + \overline{\sum_{m=1}^{\infty} b_m z^m}$ be a sense preserving K-quasiconformal harmonic mapping in $\mathbb{D},$ such that $h$ is analytic univalent in $\mathbb{D}.$ Then
\begin{equation}\label{eq:3.1}
\sum_{m=1}^\infty |a_m| r^m+\sum_{m=1}^\infty |b_m| r^m \leq 1 
\end{equation}
for $|z|=r\leq r_0=(2K+1-\sqrt{K(3K+2)})/(K+1)$ and it is sharp for $p(z)=z/(1-z)^2+k \overline{z/(1-z)^2}.$ If we take $h$ to be convex univalent then the inequality in (\ref{eq:3.1}) holds for $|z|=r\leq r_0=(K+1)/(3K+1)$ and it is sharp for $p(z)=z/(1-z)+k \overline{z/(1-z)},$ where $k=(K-1)/(K+1).$
Further, let $f_1(z)=h_1(z)+\overline{g_1(z)}=\sum_{m=1}^{\infty} c_m z^m + \overline{\sum_{m=0}^{\infty} d_m z^m} \in S(f).$ Then
\begin{equation}\label{eq:3.2}
\sum_{m=1}^\infty |c_m| r^m+\sum_{m=1}^\infty |d_m| r^m \leq 1 
\end{equation}
for $|z|=r\leq r_0=min(1/3,(2K+1-\sqrt{K(3K+2)})/(K+1)).$ If we take $h$ to be convex univalent then the inequality in (\ref{eq:3.2}) holds for $|z|=r\leq r_0=min(1/3,(K+1)/(3K+1)).$ 
\end{theorem}

In next theorem, we establish Bohr's phenomenon for univalent harmonic functions $f=h+\overline{g}\in S_H$ whose dilatation $g'/h'$ is suitably chosen. 
\begin{theorem}
Let $f(z)=h(z)+\overline{g(z)}=z+\sum_{m=2}^{\infty} a_m z^m + \overline{\sum_{m=1}^{\infty} b_m z^m}$ be a univalent and K-quasiconformal harmonic mapping in $\mathbb{D},$ where $h$ is analytic univalent in $\mathbb{D}$ and $g'(z)/h'(z)=ke^{i\theta}z^n,k=(K-1)/(K+1)\in(0,1),n\in \mathbb{N},\theta \in \mathbb{R}.$ Then
 \begin{equation}
 \sum_{m=1}^\infty |a_m| r^m+\sum_{m=1}^\infty |b_m| r^m \leq 1
 \end{equation}
 for $|z|=r \leq r_0,$ where $r_0$ is the only root of the equation
 \begin{equation} \label{eq:8}
 \frac{(k+1)r}{(1-r)^2}-\frac{2nkr}{(1-r)}-kn^2log(1-r)=1
 \end{equation}
 in $(0,1)$ and this $r_0$ is best possible one.
\end{theorem}
Letting $k\rightarrow 1$ (equivalently, $K\rightarrow \infty$) we obtain the following result.

\begin{corollary}
Let $f(z)=h(z)+\overline{g(z)}=z+\sum_{m=2}^{\infty} a_m z^m + \overline{\sum_{m=1}^{\infty} b_m z^m}$ be a univalent harmonic mapping in $\mathbb{D},$ where $h$ is analytic univalent in $\mathbb{D}$ and $g'(z)/h'(z)=e^{i\theta}z^n,n\in \mathbb{N},\theta \in \mathbb{R}.$ Then   
 \begin{equation}
 \sum_{m=1}^\infty |a_m| r^m+\sum_{m=1}^\infty |b_m| r^m \leq 1,a_1=1
 \end{equation}
 for $|z|=r \leq r_0,$ where $r_0$ is the only root in $(0,1)$ of the equation $\phi(r)=0,$ where
 \begin{equation}
 \phi(r)= \frac{2r}{(1-r)^2}-\frac{2nr}{(1-r)}-n^2log(1-r)-1.
 \end{equation}
This $r_0$ is the best possible one.
\end{corollary}
\begin{figure}
\centering
\includegraphics[scale=0.7]{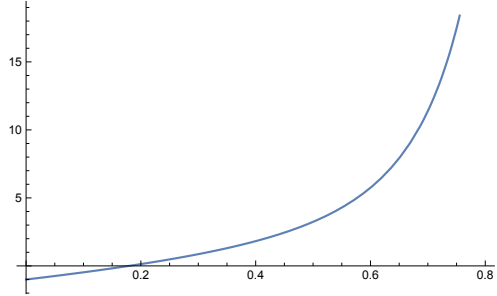}

\caption{$\phi(r)$ w.r.t $r$ for $n=3.$}

\label{fig:1}

\end{figure}

By plotting the graph of $\phi(r)$ w.r.t $r$ for different values of $n,$ we observe that there is only one root of $\phi(r)$ in $(0,1)$ which is the Bohr radius for that value of $n$ in the dilatation function. Figure \ref{fig:1} illustrates the case when $n=3$ and in the following table we have listed values of $r_0$ computed for $n=1,2,3$ and $4.$
\begin{center}
\begin{tabular}{| c  | c |}
\hline  
\hspace{1cm} $n$ \hspace{1cm} & \hspace{1cm} $r_0$ \hspace{1cm} \\
\hline
1 & 0.3485... \\
\hline
2 & 0.3121... \\
\hline
3 & 0.1794... \\
\hline
4 & 0.0959...  \\
\hline
\end{tabular}
\end{center}
We observe that if $n\rightarrow \infty,$ then $r_0\rightarrow 0.$\\
Lemma 2.1 and Theorem 2.4 together lead us to the following result for the subordination class $S(f).$
\begin{corollary}
Let $f_1(z)=h_1(z)+\overline{g_1(z)}=\sum_{m=1}^{\infty} c_m z^m + \overline{\sum_{m=1}^{\infty} d_m z^m} \in S(f)$ where $f$ is as defined in Theorem 2.4. Then
\begin{equation}
 \sum_{m=1}^\infty |c_m| r^m+\sum_{m=1}^\infty |d_m| r^m \leq 1
 \end{equation}
 for $|z|=r \leq r_1=min( 1/3, r_0), $ where $r_0$ is same as obtained in Theorem 2.4.
\end{corollary}
Next theorem shows the existence of Bohr's phenomenon for $f\in S_H$ with dilatation $w_f=(a+z)/(1+az),a\in(-1,1).$
 
\begin{theorem}
  Let $f(z)=h(z)+\overline{g(z)}=z+\sum_{m=2}^{\infty} a_m z^m + \overline{\sum_{m=1}^{\infty} b_m z^m}$ be a univalent harmonic mapping in $\mathbb{D},$ where $h$ is univalent in $\mathbb{D}$ and $g'(z)/h'(z)=\frac{a+z}{1+ az},a\in(-1,1)$ Then
 \begin{equation}
 \sum_{m=1}^\infty |a_m| r^m+\sum_{m=1}^\infty |b_m| r^m \leq 1+|a|
 \end{equation}
 for $|z|=r \leq r_0,$ where $r_0=0.2291...$ is a unique root lying in $(0,1)$ of $r^3-3r^2+5r-1= 0.$
\end{theorem} 
\begin{remark}
We observe that if we take $g'/h'=\frac{a-z}{1-az},a \in (-1,1),$ in Theorem 2.7, then we obtain the same value of $r_0.$
\end{remark}
In the following theorem we establish Bohr's phenomenon for univalent harmonic functions convex in one direction.
\begin{theorem}
Let $f(z)=h(z)+\overline{g(z)}=z+\sum_{m=2}^{\infty} a_m z^m + \overline{\sum_{m=1}^{\infty} b_m z^m}$ be a harmonic mapping in $\mathbb{D},$ where $h$ is analytic univalent in $\mathbb{D}$ and $h(z)+e^{i\theta}g(z)$ is convex univalent in $\mathbb{D}$ for some $\theta\in \mathbb{R}.$ Then
 \begin{equation}
 \sum_{m=1}^\infty |a_m| r^m+\sum_{m=1}^\infty |b_m| r^m \leq 1
 \end{equation}
 for $|z|=r \leq r_0= 0.2192....$ 

\end{theorem}

We can drop the condition of univalency of $h$ in Theorem 2.9 if we take $b_1=0.$ 

\begin{theorem}
Let $f(z)=h(z)+\overline{g(z)}=z+\sum_{m=2}^{\infty} a_m z^m + \overline{\sum_{m=2}^{\infty} b_m z^m}\in S_H^0$ be a harmonic mapping in $\mathbb{D},$ where $h(z)+e^{i\theta}g(z)$ is convex univalent in $\mathbb{D}$ for some $\theta\in \mathbb{R}.$ Then
 \begin{equation} \label{eq:2.6}
 \sum_{m=1}^\infty |a_m| r^m+\sum_{m=2}^\infty |b_m| r^m \leq 1
 \end{equation}
 for $|z|=r \leq r_0=0.3134...,$ where $r_0$ is a unique root in $(0,1)$ of $4r^3-9r^2+12r-3= 0.$ This result is sharp for Koebe function $K(z)=\frac{z-1/2z^2+1/6z^3}{(1-z)^3}+\overline{\frac{1/2z^2+1/6z^3}{(1-z)^3}}.$

\end{theorem}
Our last theorem gives Bohr radius for convex univalent harmonic functions in $S_H^0.$

\begin{theorem}
Let $f(z)=h(z)+\overline{g(z)}=z+\sum_{m=2}^{\infty} a_m z^m + \overline{\sum_{m=2}^{\infty} b_m z^m}\in K_H^0, z\in \mathbb{D}.$ Then
 \begin{equation}
 \sum_{m=1}^\infty |a_m| r^m+\sum_{m=2}^\infty |b_m| r^m \leq 1
 \end{equation}
 for $|z|=r \leq r_0=(3-\sqrt{5})/2= 0.3819.... $ This value of $r_0$ is sharp for $L(z)=\frac{1}{2}\left[\frac{z}{1-z}+\frac{z}{(1-z)^2} +\overline{\frac{z}{1-z}-\frac{z}{(1-z)^2}}\right].$
 \end{theorem}
 \section{Proof of theorems}
 We begin this section by stating a lemma which is easy to prove.
 \begin{lemma}
Let $h(z)=\sum_{m=0}^{\infty} a_m z^m$ and $g(z)=\sum_{m=0}^\infty b_m z^m$ be two holomorphic functions in $\mathbb{D}$ such that $h(z)=g(z).$ Then 
\begin{equation} \label{eq:10}
\sum_{m=0}^\infty |a_m| r^m = \sum_{m=0}^{\infty} |b_m| r^m
\end{equation}
for all $|z|=r<1.$   
\end{lemma}
\textbf{Proof of Theorem 2.4} From $g'(z)=ke^{i\theta} z^n h'(z),$ we get
$$\sum_{m=1}^\infty m b_m z^{m-1}= ke^{i\theta}\sum_{m=1}^{\infty} m a_m z^{n+m-1}, z\in \mathbb{D},$$
where $a_1=1$ and on integrating we obtain
$$ \sum_{m=1}^\infty b_m z^{m}= ke^{i\theta}\sum_{m=1}^{\infty} \frac{m}{m+n} a_m z^{m+n}, z \in \mathbb{D}.$$
Now, applying Lemma 3.1,  we get 
\begin{equation} \label{eq:11}
\sum_{m=1}^\infty |b_m| r^{m}= k\sum_{m=1}^{\infty} \frac{m}{m+n} |a_m| r^{m+n}
\end{equation} 
for all $|z|=r<1.$ Since $h$ is analytic univalent in $\mathbb{D}$ and according to De Brange's theorem, $|a_m|\leq m, m=2,3,...,$ therefore, from \eqref{eq:11}, we have 
\begin{align*}
\sum_{m=1}^\infty |a_m| r^m+\sum_{m=1}^\infty |b_m| r^m &= \sum_{m=1}^\infty |a_m| r^m+k\sum_{m=1}^{\infty} \frac{m}{m+n} |a_m| r^{m+n}\\  & \leq \sum_{m=1}^\infty m r^m+k\sum_{m=1}^{\infty} \frac{m^2}{m+n}  r^{m+n}\\ & = \sum_{m=1}^\infty m r^m+k\sum_{m=n+1}^{\infty} \frac{(m-n)^2}{m}  r^{m}\\ & \leq \sum_{m=1}^\infty m r^m+k\sum_{m=1}^{\infty} \frac{(m-n)^2}{m}  r^{m}\\ & =(k+1)\sum_{m=1}^\infty mr^m +kn^2 \sum_{m=1}^\infty\frac{1}{m}r^m-2kn\sum_{m=1}^\infty r^m \\&= \frac{(k+1)r}{(1-r)^2}-kn^2\log(1-r)-\frac{2knr}{1-r}.
\end{align*}
Thus $\sum_{m=1}^\infty |a_m| r^m+\sum_{m=1}^\infty |b_m| r^m \leq 1$ if 
\begin{equation}\label{eq:12}
\frac{(k+1)r}{(1-r)^2}-kn^2\log(1-r)-\frac{2knr}{1-r} \leq 1.
\end{equation} 
Now, we need to verify that inequality \eqref{eq:12} holds for $r\leq r_0,$ where $r_0$ is the unique root of the equation \eqref{eq:8} lying in (0,1). For this let 
$$\phi(r)=\frac{(k+1)r}{(1-r)^2}-kn^2\log(1-r)-\frac{2knr}{1-r}-1.$$ Then $\phi(r)$ is continuous in $(0,1), \phi(0)=-1<0$ and $lim_{r\rightarrow 1^-}\phi(r)>0$ implies that there is atleast one root of $\phi(r)=0$ in $(0,1).$ But $\phi'(r)>0$ for all $r\in (0,1),k\in (0,1)$ and for all $n\in \mathbb{N}$  shows that $\phi$ is strictly increasing in $(0,1).$ Hence $\phi(r)=0$ has a unique root $r_0$ in $(0,1).$\\

\begin{figure} 
\centering
\includegraphics[scale=0.8]{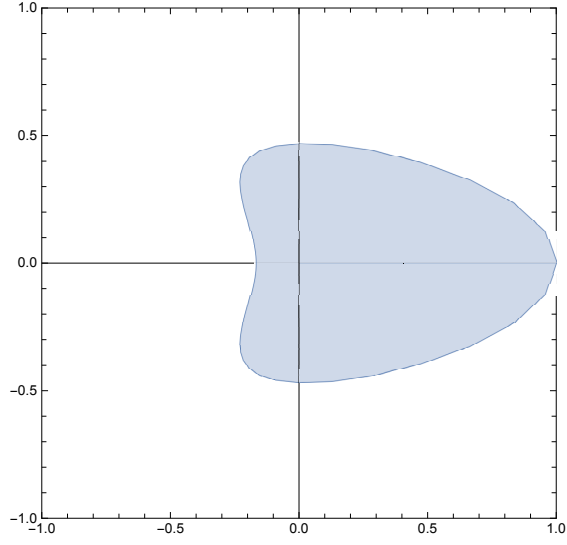}
\caption{Image of $|z|<0.3485$ under $f_0(z).$}
\label{fig:2}
\end{figure}

To see that this $r_0$ is best possible one, we consider $f_0(z)=z/(1-z)^2+\overline{z/(1-z)^2-2z/(1-z)-\log(1-z)}.$ $f_0$ maps $|z|<0.3485...$ onto the region given in the Figure \ref{fig:2} from which it is evident that $r_0$ is sharp and can not be improved further.

\vspace{0.2cm}

\noindent \textbf{Proof of Theorem 2.7} From $g'(z)=\left(\frac{a+z}{1+az}\right) h'(z)$ we obtain 
\begin{equation}
\sum_{m=1}^\infty m b_m z^{m-1}= \left(\frac{a+ z}{1+az}\right)\sum_{m=1}^{\infty} m a_m z^{m-1}, z\in \mathbb{D}
\end{equation}
where $a_1=1$ and this gives $$\sum_{m=1}^\infty m b_m z^{m-1} + \sum_{m=1}^\infty m a b_m z^{m}= \sum_{m=1}^{\infty} m a a_m z^{m-1} + \sum_{m=1}^{\infty} m a_m z^{m}, z\in \mathbb{D}.$$
Thus we have
$$\sum_{m=1}^\infty m |b_m| |z|^{m-1} - \sum_{m=1}^\infty m |a| |b_m| |z|^{m} \leq \sum_{m=1}^{\infty} m |a| |a_m| |z|^{m-1} +\sum_{m=1}^{\infty} m |a_m| |z|^{m}.$$
On integrating from 0 to r, we get
$$\sum_{m=1}^\infty  |b_m| r^{m} - \sum_{m=1}^\infty \frac{m}{m+1} |a| |b_m| r^{m+1} \leq \sum_{m=1}^{\infty}  |a| |a_m| r^{m} +\sum_{m=1}^{\infty} \frac{m}{m+1} |a_m| r^{m+1},$$
and this implies
\begin{equation} \label{eq:16}
\sum_{m=1}^\infty  \left(|b_m| - (\frac{m-1}{m})|a| |b_{m-1}|\right) r^{m} \leq \sum_{m=1}^{\infty} \left(|a||a_m|  + (\frac{m-1}{m}) |a_{m-1}|\right) r^{m}.
\end{equation}
Now, we have 
\begin{align*}
\sum_{m=1}^\infty (|a_m|+|b_m|)r^m &= \sum_{m=1}^\infty (|a_m|+|b_m|)r^m -\sum_{m=1}^{\infty} \left(\frac{m-1}{m} \right)|a||b_{m-1}|r^{m-1}+\sum_{m=1}^{\infty} \left(\frac{m-1}{m} \right)|a||b_{m-1}|r^{m-1}\\
&\leq \sum_{m=1}^\infty (|a_m|+|b_m|)r^m -\sum_{m=1}^{\infty} \left(\frac{m-1}{m} \right)|a||b_{m-1}|r^{m}+\sum_{m=1}^{\infty} \left(\frac{m-1}{m} \right)|a||b_{m-1}|r^{m-1}\\
&= \sum_{m=1}^\infty |a_m|r^m +\sum_{m=1}^{\infty} \left(|b_m|-\left(\frac{m-1}{m} \right)|a||b_{m-1}|\right)r^{m}+\sum_{m=1}^{\infty} \left(\frac{m-1}{m} \right)|a||b_{m-1}|r^{m-1}. 
\end{align*}
From \eqref{eq:16}, we get
\begin{align*}
\sum_{m=1}^\infty (|a_m|+|b_m|)r^m &\leq  \sum_{m=1}^\infty |a_m|r^m +\sum_{m=1}^{\infty} \left(|a||a_m|+\left(\frac{m-1}{m} \right)|a_{m-1}|\right)r^{m}+\sum_{m=1}^{\infty} \left(\frac{m-1}{m} \right)|a||b_{m-1}|r^{m-1} \\
&\leq \sum_{m=1}^\infty (1+|a|)|a_m|r^m +\sum_{m=1}^{\infty} \left(\frac{m-1}{m} \right)|a_{m-1}|r^{m-1}+\sum_{m=1}^{\infty} \left(\frac{m-1}{m} \right)|b_{m-1}|r^{m-1}\\
&=\sum_{m=1}^\infty (1+|a|)|a_m|r^m +\sum_{m=1}^{\infty} \left(\frac{m}{m+1} \right)\left(|a_{m}+|b_m|\right)|r^{m}.
\end{align*}
Therefore, we get 
$$\sum_{m=1}^\infty \left(\frac{1}{m+1} \right)(|a_m|+|b_m|)r^m \leq (1+|a|)\sum_{m=1}^\infty |a_m|r^m.$$
 Multiplying both sides with $r$ and then differentiating w.r.t $r,$ we get
\begin{equation}\label{eq:17}
\sum_{m=1}^\infty (|a_m|+|b_m|)r^m \leq (1+|a|)\sum_{m=1}^\infty (m+1)|a_m|r^m.
\end{equation}
As $h$ is univalent, so $|a_m|\leq m$ by De Branges's Theorem. From \eqref{eq:17}, we have
\begin{align*}
\sum_{m=1}^\infty (|a_m|+|b_m|)r^m &\leq (1+|a|)\sum_{m=1}^\infty (m+1)mr^m \\
&=(1+|a|)\left( \frac{r(1+r)}{(1-r)^3}+\frac{r}{(1-r)^2}\right) \leq (1+|a|)
\end{align*}
for $r^3-3r^2+5r-1\leq 0$ and this happens for $r\leq 0.2291....$ \\

\noindent\textbf{Proof of Theorem 2.9}
Let $h(z)+e^{i\theta}g(z)=\psi(z),$ where $\psi(z)=\sum_{m=1}^\infty C_m z^m$ is a convex univalent function. So, we have $|a_m+e^{i\theta}b_m|= |C_m|\leq 1$ for all $m\in \mathbb{N}.$ This implies $|b_m|\leq 1+|a_m|,m\in\mathbb{N}.$ We have with $|a_1|=1,$
\begin{align*}
\sum_{m=1}^\infty |a_m| r^m+\sum_{m=1}^\infty |b_m| r^m &\leq \sum_{m=1}^\infty |a_m| r^m+\sum_{m=1}^\infty (1+|a_m|) r^m\\
&=2\sum_{m=1}^\infty |a_m| r^m+\sum_{m=1}^\infty r^m .
\end{align*}
Since $h$ is univalent in $\mathbb{D},$ so by De Brange's Theorem, we have $|a_m|\leq m$ and hence we get

\begin{equation}
\sum_{m=1}^\infty |a_m| r^m+\sum_{m=1}^\infty |b_m| r^m \leq 1
\end{equation}
for $\frac{2r}{(1-r)^2}+\frac{r}{1-r}\leq 1$ i.e. for $2r^2-5r+1\geq 0.$ This is true for $r\leq r_0=\frac{5-\sqrt{17}}{4}=0.2192....$ \\
 
\begin{figure}
\label{fig:3}
\centering
\includegraphics[scale=0.8]{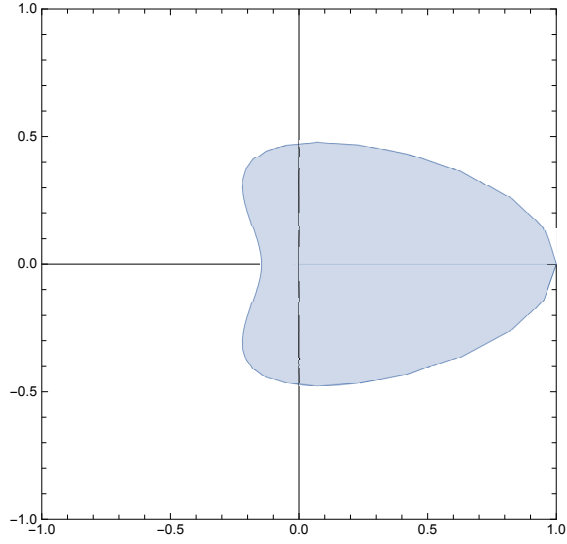}
\caption{Image of $|z|<0.3134$ under $K(z).$}
\end{figure}

 Now to prove Theorem 2.10, we first state the following result of Sheil-Small \cite{sheil1990constants}. \\
\begin{lemma}

If $f(z)=h(z)+\overline{g(z)}=z+\sum_{m=2}^{\infty} a_m z^m + \overline{\sum_{m=2}^{\infty} b_m z^m}\in S_H^0$ is convex in one direction, then 
$$|a_m|\leq \frac{(m+1)(2m+1)}{6} \hspace{2cm} |b_m|\leq \frac{(m-1)(2m-1)}{6}. $$

\end{lemma} 
\noindent\textbf{Proof of Theorem 2.10}  $h+e^{i\theta}g$ is convex univalent implies that $f$ is convex in the direction $-\theta/2,$ by the well known result of Clunie and Sheil-Small \cite{clunie1984harmonic}. Therefore, from Lemma 3.2, we have
$$|a_m|\leq \frac{(m+1)(2m+1)}{6} \hspace{2cm} |b_m|\leq \frac{(m-1)(2m-1)}{6}$$ and so, 
\begin{align*}
\sum_{m=1}^\infty |a_m| r^m+\sum_{m=2}^\infty |b_m| r^m & \leq \sum_{m=1}^\infty \frac{(m+1)(2m+1)}{6} r^m+\sum_{m=2}^\infty \frac{(m-1)(2m-1)}{6} r^m \\&=\sum_{m=1}^\infty \frac{2m^2+1}{3} r^m \\ &=\frac{2r(r+1)}{3(1-r)^3}+ \frac{r}{3(1-r)}\\ &\leq 1
\end{align*}

if $4r^3-9r^2+12r-3 \leq 0.$ This inequality holds for $r\leq r_0=0.3134...,$ where $r_0$ is unique root of $4r^3-9r^2+12r-3 = 0$ in $(0,1).$ This result is sharp for $K(z)=\frac{z-1/2z^2+1/6z^3}{(1-z)^3}+\overline{\frac{1/2z^2+1/6z^3}{(1-z)^3}},$ where $K$ is harmonic mapping in $\mathbb{D},$ which maps $|z|<0.3134$ onto region given in Figure \ref{fig:3}. It is clear from Figure \ref{fig:3} that $|K(z)|<1$ for $|z|<0.3134...$ and $0.3134...$ can not be improved. Hence this $r_0$ is sharp for inequality \eqref{eq:2.6} also.
\vspace{1mm}

\vspace{0.2cm}
To prove Theorem 2.11, we need following result of Duren \cite{duren1994harmonic}.
\begin{lemma}

If a harmonic function $f(z)=h(z)+\overline{g(z)}=z+\sum_{m=2}^{\infty} a_m z^m + \overline{\sum_{m=2}^{\infty} b_m z^m}\in K_H^0, z\in \mathbb{D},$ then 
$$|a_m|\leq \frac{m+1}{2} \hspace{2cm} |b_m|\leq \frac{m-1}{2}.$$

\end{lemma}

\noindent \textbf{Proof of Theorem 2.11} In view of Lemma 3.3 $f(z)\in K_H^0$ implies that $$|a_m|\leq \frac{m+1}{2} \hspace{2cm} |b_m|\leq \frac{m-1}{2} .$$ This gives for $a_1=1,$
\begin{align*}
\sum_{m=1}^\infty |a_m| r^m+\sum_{m=2}^\infty |b_m| r^m &\leq \sum_{m=1}^\infty \frac{m+1}{2} r^m+\sum_{m=2}^\infty \frac{m-1}{2} r^m \\
&=2\sum_{m=1}^\infty m r^m\\ &=\frac{r}{(1-r)^2} \\&\leq 1.
\end{align*}
for $r\leq r_0=0.3819...$ This value of $r_0$ is best possible, as the result is sharp for $L(z)=\frac{1}{2}\left[\frac{z}{1-z}+\frac{z}{(1-z)^2} +\overline{\frac{z}{1-z}-\frac{z}{(1-z)^2}}\right].$ 
For $L(z)$ we have
\begin{align*}
\sum_{m=1}^\infty |a_m| r^m+\sum_{m=2}^\infty |b_m| r^m &= \sum_{m=1}^\infty \left|\frac{m+1}{2}\right| r^m+\sum_{m=2}^\infty \left|\frac{1-m}{2}\right| r^m \\&=\sum_{m=1}^\infty m r^m \\&= \frac{r}{(1-r)^2} \\&\leq 1 
\end{align*}
for $r \leq 0.3819...$ Thus $r_0$ is sharp.

\medskip

\printbibliography

@article{bohr1914theorem,
  title={A theorem concerning power series},
  author={Bohr, Harald},
  journal={Proceedings of the London Mathematical Society},
  volume={2},
  number={1},
  pages={1--5},
  year={1914},
  publisher={Wiley Online Library}
}

@incollection{muhanna2017bohr,
  title={On the Bohr inequality},
  author={Muhanna, Yusuf Abu and Ali, Rosihan M and Ponnusamy, Saminathan},
  booktitle={Progress in Approximation Theory and Applicable Complex Analysis},
  pages={269--300},
  year={2017},
  publisher={Springer}
}

@article{muhanna2010bohr,
  title={Bohr's phenomenon in subordination and bounded harmonic classes},
  author={Muhanna, Yusuf Abu},
  journal={Complex Variables and Elliptic Equations},
  volume={55},
  number={11},
  pages={1071--1078},
  year={2010},
  publisher={Taylor \& Francis}
}

@article{bhowmik2019bohr,
  title={Bohr phenomenon for locally univalent functions and logarithmic power series},
  author={Bhowmik, Bappaditya and Das, Nilanjan},
  journal={Computational Methods and Function Theory},
  volume={19},
  number={4},
  pages={729--745},
  year={2019},
  publisher={Springer}
}

@article{bhowmik2018bohr,
  title={Bohr phenomenon for subordinating families of certain univalent functions},
  author={Bhowmik, Bappaditya and Das, Nilanjan},
  journal={Journal of Mathematical Analysis and Applications},
  volume={462},
  number={2},
  pages={1087--1098},
  year={2018},
  publisher={Elsevier}
}

@article{sheil1990constants,
  title={Constants for planar harmonic mappings},
  author={Sheil-Small, Terry},
  journal={Journal of the London Mathematical Society},
  volume={2},
  number={2},
  pages={237--248},
  year={1990},
  publisher={Narnia}
}

@article{clunie1984harmonic,
  title={Harmonic univalent functions},
  author={Clunie, J and Sheil-Small, T},
  journal={Ann. Acad. Sci. Fenn. Ser. A. I Math.},
  volume={9},
  pages={3--25},
  year={1984}
}

@article{duren1994harmonic,
  title={Harmonic mappings in the plane},
  author={Duren, Peter},
  journal={LECTURE NOTES IN MATHEMATICS-SPRINGER VERLAG-},
  pages={408--408},
  year={1994},
  publisher={SPRINGER VERLAG KG}
}

@article{kayumov2018bohr,
  title={Bohr radius for locally univalent harmonic mappings},
  author={Kayumov, Ilgiz R and Ponnusamy, Saminathan and Shakirov, Nail},
  journal={Mathematische Nachrichten},
  volume={291},
  number={11-12},
  pages={1757--1768},
  year={2018},
  publisher={Wiley Online Library}
}

@article{schaubroeck2000subordination,
  title={Subordination of planar harmonic functions},
  author={Schaubroeck, Lisbeth E},
  journal={Complex Variables and Elliptic Equations},
  volume={41},
  number={2},
  pages={163--178},
  year={2000},
  publisher={Taylor \& Francis}
}

@article{liu2020bohr,
  title={Bohr’s phenomenon for the classes of Quasi-subordination and K-quasiregular harmonic mappings},
  author={Liu, Ming-Sheng and Ponnusamy, Saminathan and Wang, Jun},
  journal={Revista de la Real Academia de Ciencias Exactas, F{\'\i}sicas y Naturales. Serie A. Matem{\'a}ticas},
  volume={114},
  number={3},
  pages={1--15},
  year={2020},
  publisher={Springer}
}

@article{liu2019bohr,
  title={Bohr radius for subordination and K-quasiconformal harmonic mappings},
  author={Liu, ZhiHong and Ponnusamy, Saminathan},
  journal={Bulletin of the Malaysian Mathematical Sciences Society},
  volume={42},
  number={5},
  pages={2151--2168},
  year={2019},
  publisher={Springer}
}
\end{document}